\documentclass[leqno]{article}

\title{
\textbf{Pontryagin forms on $(4r-2)$-manifolds
and symplectic structures on the spaces
of Riemannian metrics}}

\author{\textsc{R. Ferreiro P\'erez}\\
Departamento de Econom\'{\i}a
Financiera y Contabilidad I\\
Facultad de Ciencias Econ\'omicas y Empresariales\\
Universidad Complutense de Madrid\\
Campus de Somosaguas, 28223-Pozuelo de Alarc\'on,
Spain\\
\emph{E-mail:} \texttt{roferreiro@ccee.ucm.es}\\
\smallskip
\and
\textsc{J. Mu\~{n}oz Masqu\'e}
Insituto de F\'{\i}sica Aplicada, CSIC\\
C/ Serrano 144, 28006-Madrid, Spain\\
\emph{E-mail:} \texttt{jaime@iec.csic.es}}

\usepackage{amsmath,amsthm,verbatim}
\usepackage[psamsfonts]{amssymb}

\newtheorem{theorem}{Theorem}[section]
\newtheorem{proposition}[theorem]{Proposition}
\newtheorem{lemma}[theorem]{Lemma}
\newtheorem{corollary}[theorem]{Corollary}

\theoremstyle{remark}
\newtheorem{remark}[theorem]{Remark}

\date{}

\begin{document}

\maketitle

\begin{abstract}
\noindent The Pontryagin forms
on the $1$-jet bundle of Riemannian metrics,
are shown to provide in a natural way
diffeomorphism-invariant pre-symplectic
structures on the space of Riemannian metrics
for the dimensions $n\equiv 2\pmod 4$.
The equivariant Pontryagin forms provide canonical
moment maps for these structures. In dimension two,
the symplectic reduction corresponding
to the pre-symplectic form and its moment map
attached to the first Pontryagin form, is proved
to coincide with the Teichm\"{u}ller space endowed
with the Weil-Petersson symplectic form.
\end{abstract}

\bigskip


\medskip


\medskip

\noindent\emph{Acknowledgments:\/} Supported by
``Ministerio de Ciencia e Innovaci\'on'' of Spain
under grant \#MTM2008--01386.

\section{Introduction}
The aim of this paper is to show how the classical
Chern-Weil construction of Pontryagin classes,
when applied to the universal Levi-Civita connection
on the principal bundle of linear frames
over the $1$-jet bundle of Riemannian metrics,
naturally provides diffeomorphism-invariant
pre-symplectic structures on the space of Riemannian
metrics in dimensions $n=4r-2$. Similarly,
the Berline-Vergne construction of equivariant
chacteristic classes, when applied to the previous
bundle, provides a canonical moment map
for these pre-symplectic structures. In dimension two,
we apply the Marsden-Weinstein symplectic reduction
to the pre-symplectic form and its moment map
corresponding to the first Pontryagin polynomial,
and we prove that the reduction coincides with
the Teichm\"{u}ller space endowed with
the Weil-Petersson symplectic form multiplied
by the scalar factor $\frac{1}{2\pi ^2}$.

The general idea is to consider the constructions
of Riemannian geometry as depending on an arbitrary
Riemannian metric in order to obtain diffeomorphisms
invariant objets. This is equivalent to work
in the product space $M\times \mathfrak{Met}M$
instead of $M$. Moreover, as usually the objects
constructed depend only on the metric and its derivatives
up to a certain order $r$, this infinite-dimensional
manifold can be replaced by the finite-dimensional bundle
of metrics $\mathcal{M}_M$, or by its $r$-jet bundle
$J^r\mathcal{M}_M$. Working on finite-dimensional
manifolds, allows us to use local differential methods
and we have the jet bundle geometry (e.g.\ see
\cite{saunders,olver}) at our disposal. The relationship
between both spaces is given by the evaluation map
$\mathrm{ev}_r\colon M\times \mathfrak{Met}M
\to J^r\mathcal{M}_M$, $\mathrm{ev}_r(x,g)=j_x^rg$.

We work on the first jet bundle $J^1\mathcal{M}_M$
due to the fact that the Levi-Civita connection
of a metric depends on the first derivatives
of that metric. Thus, if we pull the linear frame
bundle $FM$ back to $J^1\mathcal{M}_M$, then we obtain
a principal
$Gl(n,\mathbb{R})$-bundle ${\bar{\pi }}\colon q_1^\ast FM
\to J^1\mathcal{M}_M$ admitting a canonical connection
$\boldsymbol{\omega }$ (called the universal Levi-Civita
connection), which is invariant under the action of
diffeomorphisms of $M$ (for details see \cite{natconn}).
According  to the Chern-Weil theory of characteristic
classes, by evaluating the $k$-th Pontryagin polynomial
at the curvature $\boldsymbol{\Omega }$ of the universal
Levi-Civita connection, we obtain a closed
$4k$-form on $J^1\mathcal{M}_M$, the so-called universal
$k$-th Pontryagin form $p_k(\boldsymbol{\Omega })$,
which is invariant under the natural action
of the diffeomorphism group of $M$. If $4k\leq \dim M$,
then pulling $p_k(\boldsymbol{\Omega})$ back via $j^1g$,
we obtain the $k$-th Pontryagin form $p_k(\Omega ^g)$
corresponding to the metric $g$.
As the space of Riemannian metrics on $M$ is contractible,
the maps $j^1g$ for different $g$'s are homotopic. Hence,
we recover the well-known result according to which
the cohomology class of $p_k(\Omega ^g)$ is independent
of the metric $g$ chosen. Hence, for $4k\leq \dim M$,
the universal Pontryagin forms determine the Pontryagin
classes of $M$. Moreover, as $\dim J^1\mathcal{M}_M>\dim M$,
non-zero universal Pontryagin forms of degree greater
than $\dim M$, exist. These are precisely the forms
under consideration below.

According to \cite{equiconn}, a differential $r$-form
$\alpha \in \Omega ^r(J^1\mathcal{M}_M)$
on $J^1\mathcal{M}_M$ with $r>\dim M=n$ determines
a $(r-n)$-form on the space of Riemannian metrics
$\mathfrak{Met}M$ given by $\Im \lbrack \alpha]
=\int _M\mathrm{ev}_1^\ast \alpha \in \Omega ^{r-n}
(\mathfrak{Met}M)$. In particular, if $4r-n=2$, i.e.,
$\dim M=4r-2$, then $\Im [p_k(\boldsymbol{\Omega})]$
is a closed differential $2$-form on $\mathfrak{Met}M$,
i.e., a pre-symplectic structure $\sigma $
on $\mathfrak{Met}M$, which is invariant under the action
of the orientation-preserving diffeomorphism group
$\mathrm{Diff}^+M$ on $\mathfrak{Met}M$.

In addition, the pre-symplectic structure $\sigma $ admits
a canonical moment map $\mu $, which is obtained as follows.
Since the universal Levi-Civita connection is invariant
under the action of the diffeomorphism group of $M$,
the Berline-Vergne construction of equivariant characteristic
classes, provides a canonical equivariant extension
for the universal $k$-th Pontryagin form, called
the equivariant $k$-th Pontryagin form. This equivariant
extension provides an equivariant extension of $\sigma $,
which is known (e.g., see \cite{AB2}) to be equivalent
to provide a moment map $\mu$ for $\sigma $.

Therefore, for $\dim M=4r-2$, the space of Riemannian metrics
is endowed with the pre-symplectic structure $\sigma $
and the moment map $\mu $ corresponding to the universal
and equivariant $r$-th Pontryagin forms.

In \cite{anomalies} the equivariant forms obtained
in the present paper are shown to be related
to the expressions of local gravitational anomalies
given by the Atiyah-Singer index theorem for families.
Moreover, it is also shown that for the study of the problem
of locality in quantum field theory, it is essential
that these forms are obtained form forms on the jet bundle.

As we have a pre-symplectic structure and a moment map, we can
study the corresponding symplectic reduction. Below, we analyze
in detail the two-dimensional case, $\dim M=2$. Some partial
results of the analogous situation in higher dimensions
can be found in \cite{GeoD}. For a surface, we obtain concrete
expressions of $\sigma $ and $\mu $ (see Proposition \ref{exprWP}),
and the Marsden-Weinstein quotient is proved to coincide,
up to the scalar factor $\frac{1}{2\pi ^2}$,
with the Teichm\"{u}ller space of $M$ endowed
with the Weil-Petersson symplectic form.

Note that this scalar factor $\frac{1}{2\pi ^2}$ is precisely
what is needed for the cohomology class of this form
to coincide with Mumford's tautological class $\kappa _1$
(e.g., see \cite{wo}).

\section{Preliminaries on the geometry of the bundle of metrics}
\label{preliminaries}
In this section we recall some results appeared in \cite{natconn}.
Let $q\colon\mathcal{M}_M\to M$ be the bundle of Riemannian
metrics of an $n$-dimensional smooth manifold $M$, i.e.,
$\mathcal{M}_M=\{ g_x\in S^2(T_x^\ast M):g_x$ is positive definite
on $T_xM\} $, which is a convex open subset in $S^2(T^\ast M)$.
The global sections of this bundle are the Riemannian metrics on
$M$.

We denote by $\mathrm{Diff}^+M\subset \mathrm{Diff}M$ the subgroup
of orientation preserving diffeomorphisms, and we set
$\mathcal{G}=\mathrm{Diff}M\times \mathbb{R}$,
$\mathcal{G}^{+}=\mathrm{Diff}^+M\times \mathbb{R}$.
Every coordinate system $(U,x^i)$ on $M$ induces coordinates
$(q^{-1}U,x^i,y_{ij})$ on $\mathcal{M}_M$ by setting
$g_{x}=y_{ij}(g_{x})(dx^{i})_{x}\otimes(dx^{j})_{x}$,
$y_{ij}=y_{ji}$, for every Riemannian metric $g$ on $U$.
We set $(y^{ij})=(y_{ij})^{-1}$.

Let $\pi \colon FM\to M$ be the bundle of linear frames of $M$.
The lift of $\phi \in \mathrm{Diff}M$ (resp.\ $X\in \mathfrak{X}(M)$)
to $FM$ is denoted by $\tilde{\phi }\colon FM\to FM$ (resp.\
$\tilde{X} \in \mathfrak{X}(FM)$), see \cite[VI.1--2]{KN}. Let
$q_1^\ast FM=J^1\mathcal{M}_M\times _MFM$ be the pull-back of $FM$
to $J^1\mathcal{M}_M$ via $q_1\colon J^1\mathcal{M}_M\to M$.
There are two canonical projections
\[
\begin{array}
[c]{ccc}
q_1^\ast FM
&
\overset{\bar{q}_1}{\longrightarrow }
&
FM\\
{\scriptstyle \bar{\pi }}\downarrow
& &
\downarrow{\scriptstyle \pi }\\
J^1\mathcal{M}_M
&
\overset{q_1}{\longrightarrow }
&
M
\end{array}
\]
The first projection
$\bar{\pi }\colon q_1^\ast FM\to J^1\mathcal{M}_M$
is a principal $Gl(n,\mathbb{R})$-bundle
and the second projection
$\bar{q}_1\colon q_1^\ast FM\to FM$
is $Gl(n,\mathbb{R})$-equivariant.

The group $\mathrm{Diff}M$ acts naturally
on $\mathcal{M}_M$ and on $FM$; hence
it also acts on $q_1^\ast FM$.
If $\phi \in \mathrm{Diff}M$ (resp.\
$X\in \mathfrak{X}(M)$), its lift to
$q_1^\ast FM$ is
$\hat{\phi }=(\bar{\phi }^{(1)},\tilde{\phi })$
(resp.\ $\hat{X}=(\bar{X}^{(1)},\tilde{X})$), where
$\bar{\phi }\colon \mathcal{M}_M\to \mathcal{M}_M$,
$\bar{\phi }^{(1)}\colon J^1\mathcal{M}_M\to J^1\mathcal{M}_M$
(resp.\ $\bar{X}\in \mathfrak{X}(\mathcal{M}_M)$,
$\bar{X}^{(1)}\in \mathfrak{X}(J^1\mathcal{M}_M)$)
are the natural lifts of $\phi $ (resp.\ $X$).
As $\bar{q}_1$ is a $\mathrm{Diff}M$-equivariant
map, we have $(\bar{q}_1)_\ast (\hat{X})=\tilde{X}$.
If $X=X^i\partial/\partial x^i$, then we have
\begin{equation}
\bar{X}=X^i\frac{\partial }{\partial x^i}
-\sum _{i\leq j}
\left(
\frac{\partial X^r}{\partial x^i}y_{kj}
+\frac{\partial X^k}{\partial x^j}y_{ki}
\right)
\frac{\partial}{\partial y_{ij}}.
\label{expXbarra}
\end{equation}

For every $t\in \mathbb{R}$, we define
$\varphi _t\in \mathrm{Aut}FM$ (resp.\
$\bar{\varphi }_t
\in \mathrm{Diff}(J^1\mathcal{M}_M)$,
resp.\ $\hat{\varphi }_t\in \mathrm{Aut}
\left(
q_1^\ast FM
\right) $) as follows:
\begin{align*}
\varphi _{t}(u)
&
=\exp (-\tfrac{t}{2})\cdot u,\\
\bar{\varphi }_t(j_x^1g)
&
=j_x^1
\left(
\exp (t)\cdot g
\right) ,\\
\hat{\varphi }_t(j_x^1g,u)
&
=\left(
j_x^1
\left(
\exp (t)\cdot g
\right) ,
\exp (-\tfrac{t}{2})\cdot u
\right) .
\end{align*}
We denote by $\xi \in \mathfrak{X}(FM)$
(resp.\ $\bar{\xi }\in \mathfrak{X}(J^1\mathcal{M}_M)$,
resp.\ $\hat{\xi }\in \mathfrak{X}(q_1^\ast FM)$)
the infinitesimal generator of the $1$-parameter group
$(\varphi _t)$ (resp.\ $(\bar{\varphi }_t)$, resp.\
$(\hat{\varphi }_t)$) defined above. We have
$q_{1\ast }(\bar{\xi })=0$,
$\bar{q}_{1\ast }(\hat{\xi })=\xi $,
and $\bar{\pi }_\ast (\hat{\xi })=\bar{\xi }$
and
\begin{equation}
\bar{\xi }=y_{ij} \frac{\partial}{\partial y_{ij}}
+y_{ij,k}\frac{\partial }{\partial y_{ij,k}}, \label{xi}
\end{equation}
where $(x^h,y_{ij},y_{ij,k})$ is the coordinate system
induced on $J^1\mathcal{M}_M$.

The group $\mathcal{G}=\mathrm{Diff}M\times \mathbb{R}$
acts by automorphisms of the principal $Gl(n,\mathbb{R})$-bundle
${\bar{\pi }}\colon q_1^\ast FM\to J^1\mathcal{M}_M$, inducing
a $\mathcal{G}$-action on the associated bundles to $q_1^\ast FM$
(such as $q_1^\ast TM$, $q_1^\ast T^\ast M$, etc.), as well as
in the space of sections and differential forms with values
on these bundles.

The bundle $q_1^\ast TM\to J^1\mathcal{M}_M$
is endowed with a universal metric given by $\mathbf{g}
\left(
\left(
j_x^1g,X
\right) ,
\left(
j_x^1g,Y
\right)
\right)
=g_x(X,Y)$,
$\forall X,Y\in T_xM$, which is invariant under the action
of the group $\mathcal{G}=\mathrm{Diff}M\times \mathbb{R}$
defined above (cf.\ \cite{natconn}).

We denote by $\omega ^g$ the Levi-Civita connection form of $g$
and by $\nabla ^g$ the covariant derivation law on the associated
vector bundles. The $\mathfrak{gl}(n,\mathbb{R})$-valued $1$-form
on $q_1^\ast FM$ defined by
$\boldsymbol{\omega}_\mathrm{hor}(X)=\omega ^g((\bar{q}_1)_\ast
X)$, $\forall X\in T_{(j_x^1g,u)}(q_1^\ast FM)$, is a
$\mathcal{G}$-invariant connection form on the principal
$Gl(n,\mathbb{R})$-bundle $\bar{\pi }\colon q_1^\ast FM\to
J^1\mathcal{M}_M$ (see \cite{natconn}), but unfortunately, it is
not $\mathbf{g}$-Riemannian. In fact, the connection
$\boldsymbol{\omega}_\mathrm{hor}$ induces a derivation law
$\nabla ^{\boldsymbol{\omega }_\mathrm{hor}}$ on the associated
bundles to $q_1^\ast FM$, and we have $\nabla ^{\boldsymbol{\omega
}_\mathrm{hor}}\mathbf{g}=\theta $, where $\theta
=(dy_{ij}-y_{ij,k}dx^k)\otimes dx^{i}\otimes dx^j$ is the
$V(q)$-valued $1$-form determining the contact structure on
$J^1\mathcal{M}_M$ and we have used the natural identification
$V(q)\cong q^\ast S^2T^\ast M=\mathcal{M}_M\times _MS^2T^\ast M$.

Fortunately, if
$\vartheta \in \Omega ^1(J^1\mathcal{M}_M,\mathrm{End}TM)$
is the form given by
\begin{equation}
\vartheta
=\mathbf{g}^{-1}\theta =y^{aj}
\left(
dy_{ia}-y_{ia,k}dx^k
\right)
\otimes dx^i\otimes
\frac{\partial }{\partial x^j},
\label{exptheta}
\end{equation}
then we can define a connection form on $q_1^\ast FM$---called the
`universal Levi-Civita connection'---as follows:
$\boldsymbol{\omega}
 =\boldsymbol{\omega}_\mathrm{hor}
+\tfrac{1}{2}\vartheta $, which is $\mathcal{G}$-invariant and
$\mathbf{g}$-Riemannian; i.e., $\nabla ^{\boldsymbol{\omega
}}\mathbf{g}=0$ (see \cite{natconn}). Then, the connection form
$\boldsymbol{\omega }$ is reducible to a connection on the
principal $O(n)$-sub-bundle
\[
OM
=\left\{
(j_x^1g,u_x)\in q_1^\ast FM\colon u_x\,
\text{is }g_x\text{-orthonormal }
\right\}
\subset q_1^\ast FM.
\]
In fact, it is the only $\mathcal{G}$-invariant connection on $OM$
(see \cite{natconn}). We consider the usual identification (e.g.,
see \cite[II, Example 5.2]{KN}) between differential forms on the
base manifold of a principal bundle taking values in the adjoint
bundle and differential forms of the adjoint type on that
principal bundle. With this identification we have
$\boldsymbol{\Omega}, \boldsymbol{\Omega}_{\mathrm{hor}} \in
\Omega ^2 (J^1\mathcal{M}_M,\mathrm{End}TM)$, where
$\boldsymbol{\Omega}, \boldsymbol{\Omega}_\mathrm{hor}$ are the
curvature forms of $\boldsymbol{\omega}$,
$\boldsymbol{\omega}_\mathrm{hor}$, respectively.

As $V(q)\cong q^\ast S^2T^\ast M$, an element
$h\in \Omega ^0(M,S^2T^\ast M)$ determines
a vertical vector field $H$
on $\mathcal{M}_M$, defined by $H(g_x)=(g_x,h_x)$,
$g_x\in \mathcal{M}_M$. Locally,
$H=h_{ij}\partial /\partial y_{ij}$
if $h=h_{ij}dx^i\otimes dx^j$.

Some important formulas used below, are the following
(see \cite{natconn}):
\begin{align}\boldsymbol{\Omega}_\mathrm{hor}
& =\left(
d\boldsymbol{\Gamma }_{jk}^i\wedge dx^k
+\boldsymbol{\Gamma }_{as}^i
\boldsymbol{\Gamma }_{jr}^adx^s\wedge dx^r
\right)
dx^j\otimes \frac{\partial }{\partial x^i},
\label{Omegahor}\\
\boldsymbol{\Gamma }_{jk}^i
& =\tfrac{1}{2}y^{ia}
(y_{aj,k}+y_{ak,j}-y_{jk,a}),
\label{crist}\\
\boldsymbol{\Omega }
& =\left(
\boldsymbol{\Omega }_{\mathrm{hor}}
\right) _\text{\textsc{A}}
-\tfrac{1}{2}\vartheta \wedge \vartheta.
\label{Omegan}
\end{align}

The universal $k$-th Pontryagin form of $M$,
$p_k(\boldsymbol{\Omega})\in \Omega ^{4k}(J^1\mathcal{M}_M)$, is
defined as the form obtained by means of the Chern-Weil theory of
characteristic classes by applying the $k$-th Pontryagin
polynomial to the curvature $\boldsymbol{\Omega}$ of the universal
Levi-Civita connection $\boldsymbol{\omega}$. These forms are
closed, $\mathcal{G}$-invariant and satisfy the following
universal property (see \cite{natconn}): for every Riemmanian
metric $g$ we have $(j^1g)^\ast (p_k(\boldsymbol{\Omega}))
=p_k(\Omega ^g)$, where $\Omega ^g\in \Omega ^2(M,\mathrm{End}TM)$
is the curvature form of the Levi-Civita connection of the metric
$g$. Hence the Pontryagin forms of degree equal to or less than
$n$ determine the Pontryagin classes of $M$. However, the key
point is that there are non-zero Pontryagin forms of degree
greater than $n$ (as $\dim (J^1\mathcal{M}_M)>n$). For example, as
$p_1(X)=-\tfrac{1}{8\pi ^2}\mathrm{tr}(X^2)$, $X\in
\mathfrak{so}(n)$, we have
\[
p_1(\boldsymbol{\Omega })
=-\tfrac{1}{8\pi ^2}\mathrm{tr}
\left(
\boldsymbol{\Omega}\wedge \boldsymbol{\Omega }
\right)
\in \Omega ^4(J^1\mathcal{M}_M).
\]
For $n=2$ this form does not vanish, as we see below.

We denote the covariant differential
of $X\in \mathfrak{X}(M)$ with respect
to the Levi-Civita connection of a metric
$g$ on $M$ by
$\nabla ^gX\in \Omega ^1(M,TM)\cong \Omega ^0
(M,\mathrm{End}TM)$.

Similarly, the covariant differential of a tensor field
$h\in \Omega ^0(M,\bigotimes ^rT^\ast M)$ is denoted by
$\nabla ^gh\in \Omega ^1
(M,{\textstyle\bigotimes ^r} T^\ast M)
\cong \Omega ^0(M,\bigotimes ^{r+1}T^\ast M)$.

Let $\eta ^g\colon {\textstyle\bigotimes ^3}
T^\ast M\to T^\ast M\otimes\mathrm{End}TM$
be the vector-bundle homomorphism $\eta ^g
(a\otimes b\otimes c) =b\otimes a\otimes c^\sharp $.
We set
$\dot{\nabla }^gh
=\eta ^g(\nabla ^gh)\in \Omega ^1(M,\mathrm{End}TM)$.

We denote by $(\nabla ^gX)_{\text{\textsc{S}}}$ and
$(\nabla ^gX)_{\text{\textsc{A}}}$ the $g$-symmetric
and $g$-skew-symmetric parts of $\nabla ^gX$ respectively,
and a similar notation is also used for $\dot{\nabla }^gh$.
We obtain $\delta ^gh=\mathrm{tr}(\dot{\nabla }^gh)
\in \Omega ^1(M)$, where $\delta ^gh$ is the $g$-divergence
of $h$. If $X=X^i\partial/\partial x^i$,
$h=h_{ij}dx^i\otimes dx^j$, then

\begin{align*}
\nabla ^gX
& =\left( \frac{\partial X^i}{\partial x^j}
+\Gamma _{jk}^iX^k
\right)
dx^j\otimes
\frac{\partial }{\partial x^i},\\
\dot{\nabla }^gh
& =g^{jb}
\left(
\frac{\partial h_{ib}}{\partial x^k}
-h_{ab}\Gamma _{ik}^a-h_{ia}\Gamma _{bk}^a
\right)
dx^{i}\otimes
\left(
dx^k\otimes \frac{\partial }{\partial x^j}
\right) ,\\
\delta ^gh
& =g^{jb}
\left(
\frac{\partial h_{ib}}{\partial x^j}
-h_{ab}\Gamma _{ij}^a-h_{ia}\Gamma _{bj}^a
\right)
dx^i.
\end{align*}

For every $X\in \mathfrak{X}(M)$ we denote by
$q_1^\ast X\in \Omega ^0(J^1\mathcal{M}_M,TM)$
the section given by
$(q_1^\ast X)(j_x^1g)=(j_x^1g,X_x)$.

The covariant differential
$\nabla ^{\boldsymbol{\omega }_{\mathrm{hor}}}
(q_1^\ast X)\in \Omega ^1(J^1\mathcal{M}_M,TM)$
is $q_1$-horizontal; hence, it can be viewed
as a section $\boldsymbol{\nabla }X$ of the bundle
$q_1^\ast \mathrm{End}TM$. Similarly, if
$h\in \Omega ^0(M,{\textstyle\bigotimes ^2}T^\ast M)$,
then
$\nabla ^{\boldsymbol{\omega }_{\mathrm{hor}}}
(q_1^\ast h)\in \Omega ^1
(J^1\mathcal{M}_M,{\textstyle\bigotimes ^2}T^\ast M)$
is $q_1$-horizontal, and
$\boldsymbol{\dot{\nabla }}h=\eta ^\mathbf{g}
\left(
\nabla ^{\boldsymbol{\omega }_{\mathrm{hor}}}
(q_1^\ast h)
\right) $
can be considered to be an element
of $\Omega ^0 (J^1\mathcal{M}_M,T^\ast M
\otimes \mathrm{End}TM) \subset \Omega ^1
(J^1\mathcal{M}_M,\mathrm{End}TM)$.

The following Lemma will be used to obtain
the explicit expression of the pre-symplectic forms.
\begin{lemma}
\label{iHO}
If $H$ is the vertical vector field on $\mathcal{M}_M$
determined by the section $h$ of $S^2T^\ast M$, then
$\iota _{H^{(1)}}\boldsymbol{\Omega }
=(\boldsymbol{\dot{\nabla }}h)_{\text{\textsc{A}}}
-(\mathbf{g}^{-1}
\left( q_1^\ast h
\right)
\circ\vartheta)_{\text{\textsc{A}}}$.
\end{lemma}

\begin{proof}
Given $j_x^1g\in J^1\mathcal{M}_M$, we consider
a normal system of coordinates for $g$ at $x$.
By virtue of \eqref{Omegan} and \eqref{Omegahor},
we have
\begin{align*}
\boldsymbol{\Omega }_{j_x^1g}
\! & =\! \tfrac{1}{2}
\left(
d\boldsymbol{\Gamma }_{jk}^i \! \wedge \! dx^k\!
-\!d\boldsymbol{\Gamma }_{ik}^j\! \wedge \! dx^k
\! -\! \left(
dy_{ia}\!\wedge\!dy_{aj}
\right) \!
\right) _{j_x^1g}
\! \otimes \!
\left(
\! dx^i\! \otimes \!
\frac{\partial }{\partial x^j}\!
\right) _{j_x^1g}\\
\! & =\! \tfrac{1}{2}
\left(
dy_{ki,j}\! \wedge \! dx^k \!
-\! dy_{kj,i} \! \wedge \! dx^k\!
-\! \left(
dy_{ia}\!\wedge \! dy_{aj}
\right) \!
\right) _{j_x^1g} \! \otimes \!
\left(
\! dx^i\! \otimes \!
\frac{\partial}{\partial x^j} \!
\right) _{j_x^1g}.
\end{align*}
By contracting this form with
$H^{(1)}=h_{ij}\partial/\partial y_{ij}
+(\partial h_{ij}/\partial x^k)
\partial/\partial y_{ij,k}$,
we obtain
\begin{align*}
\left(
\iota _{H^{(1)}}\boldsymbol{\Omega }
\right) _{j_x^1g}
\!\! & =\!\!
\tfrac{1}{2}
\left( \!
\left( \!
\frac{\partial h_{ki}}{\partial x^j}
\! -\! \frac{\partial h_{kj}}{\partial x^i}
\!\right)
dx^k\! -\! h_{ia}dy_{aj}\! + \!h_{aj}dy_{ia}
\! \right) _{j_x^1g}
\!\! \otimes \!\!
\left(
\!dx^i\! \otimes \!
\frac{\partial }{\partial x^j}
\! \right) _{j_x^1g}\\
\!\! & =\!\!
\left(
(\boldsymbol{\dot{\nabla}}h)_{\text{\textsc{A}}}
-\left(
\mathbf{g}^{-1}h\circ\vartheta
\right) _{\text{\textsc{A}}}
\right) _{j_x^1g}.
\end{align*}
\end{proof}

\section{Pre-symplectic structures on $\mathfrak{Met}M$}

Let $p\colon E\to M$ be a locally trivial fibre bundle
over a compact connected and oriented $n$-manifold
without boundary. In \cite{equiconn}, a map
$\Im \colon \Omega ^{n+k}(J^{r}E)\to \Omega ^k(\Gamma (E))$
has been defined, which provides a geometrical interpretation
of the forms on the jet bundle with degree greater
than the dimension of the base manifold.
If
$\mathrm{ev}_r\colon M\times \Gamma (E)\to J^rE$
is the evaluation map $\mathrm{ev}_r(x,s)=j_x^rs$,
then
$\Im\lbrack \alpha]
=\int _M\mathrm{ev}_r^\ast \alpha
\in \Omega ^k(\Gamma (E))$.
The map $\Im $ commutes with the exterior differential
and the action of the automorphims of the bundle,
and hence maps closed (resp.\ invariant) forms to closed
(resp.\ invariant) forms.

Let $\mathfrak{Met}M=\Gamma (M,\mathcal{M}_M)$
denote the space of Riemannian metrics on $M$.
As $\mathfrak{Met}M$ is an open subset of
$\mathcal{S}^2(M)=\Gamma (M,S^2T^\ast M)$,
we have the canonical identification
$T_g\mathfrak{Met}M\cong \mathcal{S}^2(M)$ for any
$g\in \mathfrak{Met}M$.
The group
$\mathcal{G}=\mathrm{Diff}M\times \mathbb{R}$
acts in a natural way on $\mathfrak{Met}M$
by setting,
\begin{align*}
\mathcal{G}\times \mathfrak{Met}M
&
\to \mathfrak{Met}M,\\
((\phi ,t),g)
&
\mapsto\exp (t)\cdot
\left(
\phi ^{-1}
\right) ^\ast g.
\end{align*}

In the rest of this section, we assume $\dim M=n=4r-2$,
$r\in \mathbb{N}$, and $f\in\mathcal{I}_{2r}^{O(n)}$
denotes a Weil polynomial of degree $2r$. Hence
$f(\boldsymbol{\Omega })\in \Omega ^{n+2}
(J^1\mathcal{M}_M)$ and $\sigma =\Im \lbrack f
(\boldsymbol{\Omega })]\in \Omega ^2(\mathfrak{Met}M)$
is a $\mathcal{G}^+$-invariant pre-symplectic form
on the space $\mathfrak{Met}M$ of Riemannian metrics.
The explicit expression for $\sigma $
is as follows:

\begin{theorem}
\label{theorem}
For every $g\in \mathfrak{Met}M$,
$h,k\in T_g\mathfrak{Met}M\cong\mathcal{S}^2(M)$
we have

\begin{multline*}
\sigma _g(h,k)=-2r(2r-1)\int _Mf
\left(
(\dot{\nabla }^gh)_{\text{\textsc{A}}},
(\dot{\nabla} ^gk)_{\text{\textsc{A}}},
\Omega ^g, \overset{(2r-2}{\ldots \ldots },\Omega ^g
\right) \\
-2r\int _Mf
\left(
\left(
g^{-1}h\circ g^{-1}k
\right) _{\text{\textsc{A}}},
\Omega ^g, \overset{(2r-1}{\ldots \ldots },
\Omega ^g
\right) .
\end{multline*}

\end{theorem}

\begin{proof}
By virtue of \cite[Proposition 11]{equiconn}
we have

\[
\sigma _g(h,k)=\int _M
\left(
j^1g
\right) ^\ast
\left(
\iota _{K^{(1)}}\iota _{H^{(1)}}f(\boldsymbol{\Omega},
\dotsc,\boldsymbol{\Omega})
\right) .
\]
Moreover, using Lemma \ref{iHO} we obtain

\begin{align*}
 \iota _{K^{(1)}}\iota _{H^{(1)}}f(\boldsymbol{\Omega},
\dotsc,\boldsymbol{\Omega})
& =2r \iota _{K^{(1)}}f
\left(
\iota _{H^{(1)}} \boldsymbol{\Omega},
\boldsymbol{\Omega},\dotsc,\boldsymbol{\Omega}
\right)
 \\
& =2r \iota _{K^{(1)}} f
\left((
\boldsymbol{\dot \nabla }h)_{\text{\textsc{A}}}
-\left(
\mathbf{g}^{-1}
\left( q_1^\ast h
\right)
\circ \vartheta
\right) _{\text{\textsc{A}}},
\boldsymbol{\Omega},\dotsc, \boldsymbol{\Omega}
\right)
 \\
& =2rf
\left(
-\left(
\mathbf{g}^{-1}
\left(
q_1^\ast h
\right)
\circ \mathbf{g}^{-1}
\left(
q_1^\ast k
\right)
\right) _{\text{\textsc{A}}},
\boldsymbol{\Omega},\dotsc,\boldsymbol{\Omega}
\right)
\\
& -2r(2r-1)f
\left(
(\boldsymbol{\dot {\nabla}}h)_{\textsc{A}}
-\left(
\mathbf{g}^{-1}
\left( q_1^\ast h
\right)
\circ \vartheta
\right) _{\textsc{A}},
\right.
\\
& \left.
\qquad
\qquad
\qquad
(\boldsymbol{\dot{\nabla }}k)_{\textsc{A}}
-\left(
\mathbf{g}^{-1}
\left(
q_1^\ast k
\right)
\circ \vartheta
\right) _{\textsc{A}},
\boldsymbol{\Omega},\dotsc, \boldsymbol{\Omega }
\right) .
\\
\end{align*}
Using $(j^1g)^\ast \vartheta=0$,
$(j^1g)^\ast \boldsymbol{\Omega }=\Omega ^g$,
and
$(j^1g)^\ast (\boldsymbol{\dot{\nabla}}h)_{\textsc{A}}
=(\dot{\nabla }^gh)_{\textsc{A}}$, the result follows.

\end{proof}

If $M$ is an oriented compact connected surface
and $f=p_1$ is the first Pontryagin polynomial,
then
\begin{align}
\sigma _g(h,k)
& =\tfrac{1}{4\pi ^2}\int _M\mathrm{tr}
\left(
\left(
g^{-1}h\circ g^{-1}k
\right) _{\text{\textsc{A}}}\wedge \Omega ^g
\right)
+\tfrac{1}{4\pi ^2}\int _M\mathrm{tr}
\left(
(\dot{\nabla} ^gh)_{\text{\textsc{A}}}
\wedge(\dot{\nabla }^gk)_{\text{\textsc{A}}}
\right). \label{symplectic2}
\end{align}
for every $g\in \mathfrak{Met}M$,
$h,k\in T_g\mathfrak{Met}M\cong \mathcal{S}^2(M)$.
(A simpler expression for this form is obtained
in Proposition \ref{exprWP} below.)

If $\dim M=6$, i.e., $r=2$, the basic Weil polynomials
of degree $4$, are $p_2$ and $(p_1)^2$. If we set
$t_k(X)=\mathrm{tr}(X^k)$,
$X\in \mathfrak{so}(6;\mathbb{R})$,
then
\begin{align*}
p_{1}
&
=-\tfrac{1}{8\pi ^2}t_2,\\
p_2
&
=\tfrac{1}{128\pi ^4}
\left[
\left(
t_2
\right) ^2-2t_4
\right] .
\end{align*}
For $f=t_4$, the formula in Theorem \ref{theorem}
yields,
\begin{align*}
\sigma _g(h,k)
&
=-12\int _M\mathrm{tr}
\left(
(\dot{\nabla }^gh)_{\text{\textsc{A}}}
\wedge(\dot{\nabla }^gk)_{\text{\textsc{A}}}
\wedge\Omega ^g\wedge \Omega ^g
\right) \\
&
-4\int _M\mathrm{tr}
\left(
\left(
g^{-1}h\cdot g^{-1}k
\right) _{\text{\textsc{A}}}
\wedge\Omega ^g\wedge \Omega ^g\wedge \Omega ^g
\right) ,
\end{align*}
and for $f=(t_2)^2$, similarly we obtain
\begin{align*}
\sigma _g(h,k)
&
=-4\int _M\mathrm{tr}
\left(
(\dot{\nabla }^gh)_{\text{\textsc{A}}}
\wedge
(\dot{\nabla }^gk)_{\text{\textsc{A}}}
\right)
\wedge\mathrm{tr}
\left(
\Omega ^g\wedge \Omega ^g
\right) \\
&
-8\int _M\mathrm{tr}
\left(
(\dot{\nabla }^gh)_{\text{\textsc{A}}}
\wedge\Omega ^g
\right)
\wedge \mathrm{tr}
\left(
(\dot{\nabla }^gk)_{\text{\textsc{A}}}
\wedge \Omega ^g
\right) \\
&
-4\int _M\mathrm{tr}
\left(
\left(
g^{-1}h\cdot g^{-1}k
\right) _{\text{\textsc{A}}}
\wedge\Omega ^g
\right)
\wedge\mathrm{tr}
\left(
\Omega ^g\wedge \Omega ^g
\right) .
\end{align*}

\section{Equivariant Pontryagin forms \& moment maps}
First, we recall the definition of equivariant cohomology
in the Cartan model (e.g.\ see \cite{BGV, GS}).
Let a connected Lie group $G$ act on a manifold $N$
and let $\mathfrak{g}\to \mathfrak{X}(N)$, $X\mapsto X_N$
be the induced Lie algebra homomorphism, $X_N$ being the
infinitesimal generator of the flow $L_{\exp (-tX)}$
and $L_g\colon N\to N$ given by $L_g(x)=g\cdot x$,
$\forall g\in \mathcal{G}$, $\forall x\in N$. Let
$\Omega _G(N)=\mathcal{P}^\bullet
(\mathfrak{g},\Omega ^\bullet (N))^G$
be the space of $G$-invariant polynomials on
$\mathfrak{g}$ with values in $\Omega ^\bullet (N)$.
We assign degree $2k+r$ to the polynomials in
$\mathcal{P}^k(\mathfrak{g},\Omega ^r(N))$.
The space of $G$-equivariant differential $q$-forms
is
\[
\Omega _G^q(N)
=\bigoplus {}_{2k+r=q}(\mathcal{P}^k(\mathfrak{g},
\Omega ^r(N)))^G.
\]
Let
$d_c\colon \Omega _G^q(N)\to \Omega _G^{q+1}(N)$,
$(d_c\alpha)(X)=d(\alpha(X))-i_{X_{N}}\alpha(X)$,
$\forall X\in \mathfrak{g}$, be the Cartan differential,
As is well known, on $\Omega_{G}^{\bullet}(N)$
we have $(d_c)^2=0$. The $G$-equivariant cohomology of $N$
(in the Cartan model) is the cohomology of the complex
$(\Omega_{G}^{q}(N),d_c)$.

Given a $G$-invariant closed form $\sigma \in \Omega ^q(M)$,
an equivariant differential form
$\sigma ^{\# }\in \Omega _G^q(M)$ is said to be
a $G$-equivariant extension of $\sigma$ if $d_c\sigma ^{\#}=0$
and $\sigma ^{\#}(0)=\sigma $. In general, there could exist
obstructions to the existence of equivariant extensions
(e.g., see \cite{Wu}) but for the universal Pontryagin forms,
the classical construction of equivariant characteristic classes
of Berline and Vergne (see \cite{BV1, BV2, BT}) really provides
canonical equivariant extensions: As the universal Levi-Civita
connection is $\mathcal{G}$-invariant, for every
$f\in \mathcal{I}_k^{O(n)}$ the $\mathcal{G}$-equivariant
characteristic form associated to $f$ and $\boldsymbol{\omega }$,
is a $\mathcal{G}$-equivariant extension
of $f(\boldsymbol{\Omega },\overset{(k}{\ldots },
\boldsymbol{\Omega})$, given by
$f(\boldsymbol{\Omega }_{\mathcal{G}})(X,t)
=f(\boldsymbol{\Omega }-\boldsymbol{\omega}(\hat{X}+t\hat{\xi }),
\overset{(k}{\ldots},\boldsymbol{\Omega } -\boldsymbol{\omega }
(\hat{X}+t\hat{\xi }))$.

As a simple computation shows, we have

\begin{lemma}
\label{lema3}
Let $\omega \in \Omega ^1(FM,\mathfrak{gl}(n,\mathbb{R}))$
be the connection form of a linear connection $\nabla $.
For every vector field $X\in \mathfrak{X}(M)$ we have

\begin{enumerate}
\item[\emph{(1)}]
The $0$-form $\omega (\tilde{X})\in \Omega ^0
(FM,\mathfrak{gl}(n,\mathbb{R}))$
is of adjoint type.

\item[\emph{(2)}]
If $\nabla$ is symmetric, the $0$-form on $M$ with values
on $\mathrm{End}TM$ corresponding to $\omega (\tilde{X})$
coincides with $\nabla X$.

\item[\emph{(3)}]
The $0$-form
$\omega (\xi )\in \Omega ^0(FM,\mathfrak{gl}(n,\mathbb{R}))$
is of adjoint type and it corresponds to
$-\frac{1}{2}\mathrm{id}_{TM}\in \Omega ^0(M,\mathrm{End}TM)$.
\end{enumerate}
\end{lemma}

\begin{proposition}
\label{propOmega}The explicit expression
for the $\mathcal{G}$-equivariant characteristic
form associated to the Weil polynomial $f$
is as follows:
\begin{align*}
f(\boldsymbol{\Omega }_{\mathcal{G}})(X,t)
& =f\left(
\boldsymbol{\Omega}
-\left(
\boldsymbol{\nabla }X
\right) _{\text{\textsc{A}}},
\overset{(k}{\ldots }, \boldsymbol{\Omega }
-\left(
\boldsymbol{\nabla}X
\right) _{\text{\textsc{A}}}
\right) \\
& =\sum _{i=1}^k(-1)^{k-i}
\tbinom{k}{i}f(\boldsymbol{\Omega},
\overset{(i}{\ldots },\boldsymbol{\Omega },
\left(
\boldsymbol{\nabla }X
\right) _{\text{\textsc{A}}},
\overset{(k-i}{\ldots },
\left(
\boldsymbol{\nabla }X
\right)
_{\text{\textsc{A}}}).
\end{align*}
\end{proposition}

\begin{proof}
Let $\boldsymbol{\omega }$ be the universal
Levi-Civita connection form on
$q_1^\ast FM$ and let $\hat{\xi }$ be
the vector field introduced in the section
\ref{preliminaries}. We have
$\boldsymbol{\omega}(\hat{\xi })=0$,
and for every $X\in \mathfrak{X}(M)$ the form
$\boldsymbol{\omega }(\hat{X})$ is
a $0$-form of adjoint type on $q_1^\ast FM$,
and the corresponding $\mathrm{End}TM$-valued
$0$-form on $J^1\mathcal{M}_M$ is
$\left(
\boldsymbol{\nabla}X
\right) _{\text{\textsc{A}}}$.
In fact, from \eqref{expXbarra}, \eqref{xi},
\eqref{exptheta}, and \eqref{crist} we have
\begin{align*}
\vartheta
\left(
\bar{X}^{(1)}
\right)
& =-y^{bj}
\left(
\left(
y_{aj}\boldsymbol{\Gamma }_{ki}^a
+y_{ai} \boldsymbol{\Gamma }_{kj}^a
\right)
X^k +\frac{\partial X^r}{\partial x^i}y_{kj}
+\frac{\partial X^k}{\partial x^j}y_{ki}
\right)
dx^i\otimes \frac{\partial }{\partial x^b}\\
& =-2\left(
\boldsymbol{\nabla }X
\right) _{\text{\textsc{S}}},\\
\vartheta
\left(
\bar{\xi }
\right)
& =y^{ia}y_{aj}dx^j\otimes
\frac{\partial }{\partial x^i}
=\mathrm{id}_{TM}.
\end{align*}
The $0$-forms
$\boldsymbol{\omega}_{\mathrm{hor}}(\hat{X}),
\boldsymbol{\omega}_{\mathrm{hor}}(\hat{\xi })$
are of adjoint type due to the fact that
$\boldsymbol{\omega}_{\mathrm{hor}}$ is
invariant under the action of
$\mathrm{Diff}M\times\mathbb{R}$. If $\alpha
\in \Omega ^0(J^1\mathcal{M}_M,\mathrm{End}TM)$
is the $0$-form taking values in $\mathrm{End}TM$
corresponding to
$\boldsymbol{\omega}_{\mathrm{hor}}(\hat{X})$,
then from the formula
$\boldsymbol{\omega }_{\mathrm{hor}}
(\hat{X})(u,j_x^1g)
=\omega _u^g((q_1)_\ast \hat{X})
=\omega _u^g(\tilde{X})$
and Lemma \ref{lema3}, we obtain
$\alpha (j_x^1g)=(j_x^1g,(\nabla ^gX)(x))
=(\boldsymbol{\nabla }X)(j_x^1g)$. Hence
$\alpha =\boldsymbol{\nabla }X$. Accordingly,
the form corresponding to
$\boldsymbol{\omega }(\hat{X})
=\boldsymbol{\omega}_{\mathrm{hor}}(\hat{X})
+\tfrac{1}{2}\vartheta (\bar{X}^{(1)})$
is
$\left(
\boldsymbol{\nabla }X
\right) _{\text{\textsc{A}}}
=\boldsymbol{\nabla }X
-\left(
\boldsymbol{\nabla }X
\right) _{\text{\textsc{S}}}$.

If
$\beta \in \Omega ^0(J^1\mathcal{M}_M,\mathrm{End}TM)$
is the form corresponding to
$\boldsymbol{\omega }_{\mathrm{hor}}(\hat{\xi })$,
then we have
$\boldsymbol{\omega}_\mathrm{hor} (\hat{\xi })(u,j_x^1g)
=\omega _u^g((q_1)_\ast \hat{\xi }) =\omega _u^g(\xi )$.
Hence
\[
\beta(j_x^1g)
=\left(
j_x^1g,-\tfrac{1}{2}\mathrm{id}_{T_xM}
\right)
=-\tfrac{1}{2}q_1^\ast \mathrm{id}(j_x^1g),
\]
that is,
$\beta=-\tfrac{1}{2}q_1^\ast \mathrm{id}_{TM}$.
Therefore, the form corresponding to
$\boldsymbol{\omega }(\hat{\xi })$ vanishes,
since
$\beta +\tfrac{1}{2}\vartheta (\bar{\xi })=0$.
\end{proof}

For example, the first equivariant Pontryagin
form is given by
\begin{equation}
p_1(\boldsymbol{\Omega }_{\mathcal{G}})(X,t)
\! =\! -\tfrac{1}{8\pi ^2}
\left(
\mathrm{tr}
\left(
\boldsymbol{\Omega} \!\wedge\!
\boldsymbol{\Omega}
\right)
\! - \!2 \mathrm{tr}
\left(
\left(
\boldsymbol{\nabla}X
\right) _{\text{\textsc{A}}} \!\circ \!
\boldsymbol{\Omega }
\right)
\! +\! \mathrm{tr}
\left(
\left(
\boldsymbol{\nabla}X
\right) _{\text{\textsc{A}}} \! \circ \!
\left( \boldsymbol{\nabla}X
\right) _{\text{\textsc{A}}}
\right)
\right) . \label{1equipont}
\end{equation}

Finally, we recall the relationship between equivariant extensions
of a pre-symplectic form and moment maps (e.g., see \cite{AB2}).
If $\omega $ is a pre-symplectic form on $N$, then an equivariant
extension of $\omega $ is given by $\omega ^{\#}=\omega +\mu $,
where $\mu \colon \mathfrak{g}\to \Omega ^0(N)$ is a $G$-invariant
linear map satisfying $i_{X_N}\omega =d(\mu (X))$, i.e., $\mu $ is
a (co-)moment map for $\omega $. Hence, to give an equivariant
extension for a pre-symplectic form is equivalent to giving a
moment map for it.

The map $\Im$ introduced above, naturally extends
to a map on the spaces of equivariant differential
forms
$\Im \colon \Omega _{\mathcal{G}^+}^{n+k}
(J^1\mathcal{M}_M)\to \Omega _{\mathcal{G}^+}^k
(\mathfrak{Met}M)$
that commutes with the Cartan differential
(see \cite{equiconn}). By applying this map
to the equivariant Pontryagin forms, we obtain
equivariant extensions of the pre-symplectic
structures on $\mathfrak{Met}M$, or equivalently,
canonical moment maps for them, given by
\begin{align*}
\mu \colon \mathfrak{X}(M)\times \mathbb{R}
& \to \Omega ^0(\mathfrak{Met}M),\\
\mu (X,t)_g & =-2r\int _Mf
\left(
\left(
\nabla ^gX
\right) _{\text{\textsc{A}}},
\Omega ^g,\overset{(2r-1}{\ldots \ldots\,},
\Omega ^g
\right) .
\end{align*}

In the two-dimensional case, for $f=p_1,$
the formula \eqref{1equipont}
yields the following expression:
\begin{align}
\mu \colon \mathfrak{X}(M)\times \mathbb{R}
&
\to \Omega ^{0}
(\mathfrak{Met}M),
\label{moment}\\
\mu(X,t)_g
&
=\tfrac{1}{4\pi ^2}
\int _M\mathrm{tr}
\left(
\left(
\nabla
 ^gX
\right) _{\text{\textsc{A}}}
\circ\Omega ^g
\right) .
\nonumber
\end{align}
Similarly, if $\dim M=6$, then for $f=t_4$
we obtain the following moment map:
\[
\mu(X,t)_g=-4\int _M\mathrm{tr}
\left(
\left(
\nabla ^gX
\right) _{\text{\textsc{A}}}
\circ\Omega ^g
\wedge \Omega ^g
\wedge \Omega ^g
\right) ,
\]
and for $f=(t_{2})^{2}$, we have
\[
\mu(X,t)_g=-4\int_M\mathrm{tr}
\left(
\left(
\nabla ^gX
\right) _{\text{\textsc{A}}}
\circ\Omega ^g
\right)
\wedge\mathrm{tr}
\left(
\Omega ^g
\wedge\Omega ^g
\right) .
\]

\section{ Symplectic reduction in dimension 2}

From now on, we assume that $M$ is a compact,
orientable surface, so that $n=2$. By applying
the preceding considerations to the first Pontryagin
polynomial $p_1$ we obtain a canonical
$\mathcal{G}$-invariant pre-symplectic structure
$\sigma $ on the space of Riemannian metrics
$\mathfrak{Met}M$ and a moment map $\mu $ for it,
given by \eqref{symplectic2} and \eqref{moment}.
In this section we apply the Marsden-Weinstein procedure
of symplectic reduction to the pre-symplectic manifold
$(\mathfrak{Met}M,\sigma )$ with respect the moment map
$\mu $. First we obtain simpler expressions of $\sigma $
and $\mu$ in the 2-dimensional case.

The following results easily follows by a direct
calculation in normal coordinates.
\begin{lemma}
We have
\begin{enumerate}

\item[\emph{a)}]
$\mathrm{Pfaff}((\dot{\nabla }^gh)_{\textsc{A}})
=-\tfrac{1}{2}\star_g(\delta ^gh-d(\mathrm{tr}_gh))$
for every $h\in \mathcal{S}^2(M)$.

\item[\emph{b)}]
$\left(
\nabla ^gX
\right) _{\textsc{A}}
=\frac{1}{2}g^{-1}(dX^{\flat})$
for every $X\in \mathfrak{X}(M)$.
\end{enumerate}
\end{lemma}
\begin{proposition}
\label{exprWP}
The expressions of $\sigma $ and $\mu $
are as follows:
\begin{align*}
\sigma _g(h,k)
& = \tfrac{1}{4\pi ^2}
\int _MS ^g\mathrm{tr}
\left( g^{-1}h\circ g^{-1}k\circ g^{-1}\mathrm{vol}_g
\right)
\mathrm{vol}_g\\
& \quad
-\tfrac{1}{8\pi^{2}}\int _M
\left(
\delta ^gh-d(\mathrm{tr}_gh)
\right)
\wedge
\left(
\delta ^gk-d(\mathrm{tr}_gk)
\right) \\
\mu _g(X,t) & =\tfrac{1}{4\pi ^2}
\int _MdS ^g\wedge X^\flat,
\end{align*}
for every $g\in \mathfrak{Met}M$, $h,k\in \mathcal{S}^2(M)\cong
T_g\mathfrak{Met}M$, $X\in \mathfrak{X}(M)$ and $t\in \mathbb{R}$.
\end{proposition}

\begin{remark}
In the previous proposition $S ^g$ denotes the scalar curvature
of the metric $g$, i.e., twice the Gauss curvature.
\end{remark}

\begin{proof}
The curvature tensor of a Riemannian metric $g$
on a surface is given by (e.g., see \cite{KN}),
\[
\Omega ^g=S^g
\left(
g^{-1}\mathrm{vol}_g
\right)
\otimes \mathrm{vol}_g\in \Omega ^2
(M,\mathrm{End}TM),
\]
and we have
\begin{align*}
\mathrm{tr}
\left(
\left( g^{-1}h\circ g^{-1}k
\right)
_{\text{\textsc{A}}} \circ \Omega ^g
\right)
& =
S ^g\mathrm{tr}
\left(
\left(
g^{-1}h\circ g^{-1}k
\right) _{\text{\textsc{A}}}
\circ g^{-1}\mathrm{vol}_g
\right)
\mathrm{vol}_g\\
& =S ^g\mathrm{tr}
\left(
g^{-1}h\circ g^{-1}k\circ
g^{-1}\mathrm{vol}_g
\right)
\mathrm{vol}_g,
\end{align*}
where the last equality is due to the fact
that the trace of the product of a symmetric
and a skew-symmetric endomorphism vanishes.
This gives the first term in the expression
of $\sigma $. For the second term we have
\begin{align*}
\mathrm{tr}
\left(
(\dot{\nabla} ^gh)_{\textsc{A}}
\wedge(\dot{\nabla} ^gk)_{\textsc{A}}
\right)
& =-2\mathrm{Pfaff}
((\dot{\nabla}^gh)_{\textsc{A}})
\wedge \mathrm{Pfaff}
((\dot{\nabla}^gk)_{\textsc{A}})
\\
&=-\tfrac{1}{2}\star _g
\left(
\delta^gh-d(\mathrm{tr}_gh)
\right)
\wedge\star _g
\left(
\delta^gk-d(\mathrm{tr}_gk)
\right)
\\
&=-\tfrac{1}{2}
\left(
\delta^gh-d(\mathrm{tr}_gh)
\right)
\wedge
\left(
\delta^gk-d(\mathrm{tr}_gk)
\right),
\end{align*}
where the last equality is due to the fact
that, in dimension $2$, the following
formula holds:
$\star _g\alpha \wedge \star_g \beta
=\alpha \wedge \beta $
for $\alpha ,\beta \in \Omega ^1(M)$.

For the expression of the moment map,
we have
\begin{align*}
\mathrm{tr}
\left(
 \left(
 \nabla ^gX
 \right) _{\textsc{A}} \circ
\Omega ^g
\right)
&=\tfrac{1}{2}S^g\mathrm{tr}
\left(
g^{-1}dX^\flat \circ g^{-1}\mathrm{vol}_g
\right)
\mathrm{vol}_g
\\
&=-S ^g\mathrm{Pfaff}(g^{-1}dX^\flat )
\mathrm{Pfaff}(g^{-1}\mathrm{vol}_g)\mathrm{vol}_g
\\
&=-S ^gdX^\flat ,
\end{align*}
and hence
\begin{align*}
\mu _g(X,t)
& =\tfrac{1}{4\pi ^2}\int _M\mathrm{tr}
\left(
\left(
\nabla ^gX
\right) _{\text{\textsc{A}}} \circ\Omega ^g
\right)
=-\tfrac{1}{4\pi ^2}\int _MS ^gdX^\flat
=\tfrac{1}{4\pi ^2}\int
_MdS ^g\wedge X^\flat .
\end{align*}
\end{proof}

\begin{remark}
The form $\sigma $ is not a symplectic form,
as it is degenerate. In fact, we have
$\sigma _g(g,h)=0$ for every
$g,h\in \mathcal{S}^2(M)$, as $\delta ^gg=0$,
$\mathrm{tr}_gg=2$, and hence
\begin{align*}
\sigma _g(g,h)
& =\tfrac{1}{4\pi ^2}\int _MS ^g\mathrm{tr}
\left(
g^{-1}g\circ g^{-1}h\circ g^{-1}\mathrm{vol}_g
\right)
\mathrm{vol}_g\\
& =\tfrac{1}{4\pi ^2}\int _MS^g\mathrm{tr}
\left(
g^{-1}h\circ
g^{-1}\mathrm{vol}_g
\right)
\mathrm{vol}_g=0.
\end{align*}
\end{remark}

\begin{corollary}
Let $g\in \mathfrak{Met}M$ be a Riemaniann metric
on an oriented compact connected surface. Then,
$\mu _g(X,t )=0$, $\forall X\in \mathfrak{X}(M)$,
$\forall t \in \mathbb{R}$, if and only if
the scalar curvature $S ^g$ of $g$, is constant.
Hence $\mu ^{-1}(0)=\mathfrak{Met}_{\mathrm{const}}M$
is the space of metrics of constant curvature.
\end{corollary}

From now on, we assume the genus of $M$ is $\gamma >1$.
By the Gauss-Bonnet theorem, the space of metrics
of constant scalar curvature $-1$ can be identified
to $\mathfrak{Met}_{-1}M
\cong \mathfrak{Met}_{\mathrm{const}}M/\mathbb{R}$.
Hence the Marsden-Weinstein quotient
\[
\mu ^{-1}(0)/(\mathrm{Diff}^+M\times \mathbb{R)}
\cong \mathfrak{Met}_{-1}M/\mathrm{Diff}^+M
\cong \mathcal{M}_\gamma
\]
is the moduli space of complex surfaces of genus
$\gamma $.

As the moduli space presents singularities
due to the fact that the action of $\mathrm{Diff}^+M$
on $\mathfrak{Met}_{-1}M$ is not free, it is customary
to replace $\mathrm{Diff}^+M$ by the connected component
of the identity
$\mathrm{Diff}^eM\subset\mathrm{Diff}M$.
The action of $\mathrm{Diff}^eM$ on
$\mathfrak{Met}_{-1}M$ is free, and the quotient space
$\mathfrak{Met}_{-1}M/\mathrm{Diff}^eM=\mathcal{T}(M)$
is the Teichm\"{u}ller space of $M$.

The restriction of $\sigma$ to $\mathfrak{Met}_{-1}M$ projects
onto a canonical pre-symplectic form $\underline{\sigma }$
on $\mathcal{T}(M)$. Below we show that $\underline{\sigma }$
basically coincides with the Weil-Petersson symplectic form.

First let us recall some results about the Teichm\"{u}ller
space and the Weil-Petersson metric. We follow the exposition
in \cite{Tromba}.

The group $\mathrm{Diff}^{e}M$ acts properly and freely
on the manifold $\mathfrak{Met}_{-1}M$, and the quotient
space
$\mathcal{T}(M)=\mathfrak{Met}_{-1}M/\mathrm{Diff}^eM$
is a differentiable manifold of dimension $6\gamma -6$
called the Teichm\"{u}ler space of $M$. For every
$g\in \mathfrak{Met}_{-1}M$, we have the identification
\[
T_{[g]}\mathcal{T}(M)\cong \mathcal{S}^2(g)^{TT}
=\left\{
h\in \mathcal{S}^2(M):\mathrm{tr}_gh=0,
\delta ^gh=0
\right\} .
\]

On the space $\mathfrak{Met}M$ there exists
a canonical Riemannian metric $\mathbf{G}$
(see \cite{gilmichor,Tromba}) given by
$\mathbf{G}_g(h,k)=\tfrac{1}{2}\int _M\mathrm{tr}
\left(
g^{-1}h\circ g^{-1}k
\right)
\mathrm{vol}_g$, for every $g\in \mathfrak{Met}M$,
and every
$h,k\in \mathcal{S}^2(M)\cong T_g\mathfrak{Met}M$,
which is invariant under the action
of the diffeomorphisms group on $\mathfrak{Met}M$.
The metric $\mathbf{G}$ induces a Riemannian metric
on $\mathcal{T}(M)$, which coincides
with the Weil-Petersson metric.

Moreover, the manifold $\mathcal{T}(M)$ is endowed
with a complex structure $\mathcal{J}$, given by
$\mathcal{J}_{[g]}(h)=-\mathrm{vol}_gg^{-1}h$.
The metric $\mathbf{G}$ is compatible with
the complex structure $\mathcal{J}$, and
$\mathcal{T}(M)$ is a K\"{a}hler manifold. Hence
$\mathcal{T}(M)$ is endowed with a canonical
symplectic structure $\sigma _\mathrm{WP}$,
called the Weil-Petersson symplectic form, given by
\[
\left(
\sigma _{\mathrm{WP}}
\right) _{[g]}(h,k)
=\mathbf{G}_g
\left(
\mathcal{J}h,k
\right) ,
\quad
\forall g\in \mathfrak{Met}_{-1},\;
\forall h,k\in
\left(
\mathcal{S}^2(M)
\right) _g^{TT}.
\]

\begin{theorem}
\label{thsigma}
We have $\underline{\sigma}
=\tfrac{1}{2\pi ^2}\sigma _{\mathrm{WP}}$. Hence,
the symplectic reduction of
$(\mathfrak{Met}M,\sigma )$ is $(\mathcal{T}(M)$,
$\tfrac{1}{2\pi ^2}\sigma _{\mathrm{WP}})$.
\end{theorem}

\begin{proof}
By virtue of Theorem \ref{exprWP},
for every $g\in \mathfrak{Met}_{-1}M$,
$h,k\in \mathcal{S}^2(g)^{TT}$
we have
\begin{align*}
\sigma _g(h,k) &
=-\tfrac{1}{4\pi ^2}\int_M\mathrm{tr}
\left(
g^{-1}h\circ g^{-1}k\circ g^{-1}\mathrm{vol}_g
\right)
\mathrm{vol}_g\\
& =-\tfrac{1}{4\pi ^2}\int _M\mathrm{tr}
\left(
g^{-1}\mathrm{vol}_g\circ g^{-1}h\circ g^{-1}k
\right)
\mathrm{vol}_g\\
& =\tfrac{1}{2\pi ^2}\mathbf{G}_g(\mathcal{J}h,k)\\
& =\tfrac{1}{2\pi ^2}
\left(
\sigma _\mathrm{WP}
\right) _g (h,k),
\end{align*}
and the result follows.
\end{proof}

\begin{remark}
The preceding result provides an alternative proof
of the fact that the Weil-Petersson metric
on $\mathcal{T}(M)$ is K\"{a}hler, as we know that
$\underline{\sigma }$ is closed by its very definition,
and accordingly $\sigma _\mathrm{WP}$ is also closed.
\end{remark}

\begin{remark}
In \cite{Donaldson1,Donaldson2,Fujiki} the Teichm\"uller
space with the Weil-Petersson symplectic form, is obtained
by a symplectic reduction from the space of complex structures.
A comparison of our constructions with this Donaldson-Fujiki's
approach can be found in \cite{GeoD}.
\end{remark}

\end{document}